\documentclass[twoside, a4paper, 10pt]{amsart}

\usepackage{blindtext}
\usepackage{amsfonts}
\usepackage{amsthm}
\usepackage{verbatim}
\usepackage{amsmath, amssymb}
\usepackage{tikz}
\usetikzlibrary{matrix, arrows}
\usepackage{listings}
\usepackage{color}
\usepackage{listings}
\usepackage[all]{xy}
\usepackage[pdftex,colorlinks,linkcolor=blue,citecolor=blue]{hyperref}
\usepackage{graphicx}
\usepackage{float}
\usepackage[margin=3cm]{geometry}
\usepackage{bigints}
\usepackage{dsfont}
\newcommand{\Z}{\mathbb{Z}}  
\newcommand{\F}{\mathbb{F}_q}

\newcommand{\I}{\mathcal{I}}
\newcommand{\J}{\mathcal{I}}

\newcommand{\eps}{\varepsilon}

\theoremstyle{plain} 
\newtheorem{thm}{Theorem}[section]
\newtheorem{lem}{Lemma}[section]

\theoremstyle{definition}

\theoremstyle{remark} 
\newtheorem*{rema}{Remark}
\newtheorem*{ack}{Acknowledgements}

\title{On almost Cap sets in three variables and the multivariable Cap set problem}

\author{Alexander Fish}
\address{School of Mathematics and Statistics, University of Sydney, Australia}
\email{alexander.fish@sydney.edu.au}

\author{Dibyendu Roy}
\address{School of Mathematics and Statistics, University of Sydney, Australia}
\email{droy5385@uni.sydney.edu.au}

\date{\today}                                           
\begin{document}
\maketitle
\raggedbottom

\begin{abstract}
In this note we prove that almost cap sets $A \subset \F^n$, i.e., the subsets of $\F^n$ that do not contain  too many arithmetic progressions of length three, satisfy that $|A| < c_q^n$ for some $c_q < q$. As a corollary we prove a multivariable analogue  of Ellenberg-Gijswijt theorem \cite{EG}.
\end{abstract}
\section{Introduction}

We denote by $\F$ the field with $q$ elements. Let $a,b,c \in \F\setminus\{0\}$ be such that $a+b+c = 0$. A set $A \subset \F^n$ is called a \textbf{cap set} (for the tuple $(a,b,c)$) if any solution $(x,y,z) \in A^3$ of the equation $ax+by+cz = 0$ is of the form $x=y=z$. In a breakthrough paper Lev-Croot-Pach \cite{CLP} showed that for $n$ large enough, any cap set\footnote{The notion of the cap set naturally extends to modules over rings.} $A \subset \Z_4^n$ for the tuple $(1,1,-2)$ satisfies that $|A| \le r_4^n$, for some $r_4 < 4$. Extending their method, Ellenberg-Gijswijt \cite{EG} showed that for any $q$ there exists $c_q < q$ such that any cap set $A \subset \F^n$ satisfies $|A| \le c_q^n$. 

In this paper we obtain analogous upper bounds on the cardinality of \textbf{almost cap sets} in $\F^n$. These sets satisfy  a much weaker structural rigidity than cap sets.  We say that a set $A \subset \F^n$ is an $(\epsilon,\delta)-$cap set for  $\epsilon,\delta > 0$ and the tuple $(a,b,c) \in (\F \setminus \{0\})^3$ with $a+b+c = 0$, if there exists $A' \subset A$ with $|A'| > \delta |A|$ and such that for every $x \in A'$ the number of pairs $(y,z) \in A$ satisfying $ax+by+cz = 0$ is less  than $|A|^{\epsilon}$. 

Our main result gives an upper bound on the maximal cardinality of an $(\epsilon, \delta)-$cap set $A \subset \F^n$. 
For any choice of $\epsilon > 0$, $\mathbf{a} = (a,b,c) \in (\F \setminus \{0\})^3$ satisfying $a+b+c = 0$, and a set $A \subset \F^n$ we denote by $A_{\mathbf{a}}^{\epsilon}$ the following set:
\[
A_{\mathbf{a}}^{\epsilon} = \{x \in A \,|\,\, \exists \text{ at least } |A|^{\epsilon} \text{ pairs } (y,z)\in A^2  \text{ with } ax + b y + c z = 0\}.
\]

\begin{thm}
\label{cap_set_thm}
There exist $\epsilon > 0$ and $c_q < q$ such that for any $\delta > 0$,  $A \subset \F^n$ with $|A| > c_q^n$ and  $\mathbf{a} = (a,b,c) \in (\F \setminus\{0\})^3$ for sufficiently large $n$  we have 
\[
|A_{\mathbf{a}}^{\epsilon}| \ge  (1-\delta) |A|.
\]
In other words, there exists $\epsilon > 0$ such that for any $\delta >0$, the $(\epsilon,\delta)-$cap sets $A \subset \F^n$ satisfy that $|A| \le c_q^n$ for sufficiently large $n$. 
\end{thm}
\medskip
To prove Theorem \ref{cap_set_thm} we use Tao's symmetric reformulation \cite{Tao} of the method of Croot-Lev-Pach-Ellenberg-Gijswijt together with Lovett's lower bound  on the slice rank of tensors in terms of the cardinality of an independent set \cite{Lov}. 
\medskip

We also extend the Ellenberg-Gijswijt \cite{EG} upper bound on the cardinality of cap sets to the case of more than three variables.
 
\begin{thm}
\label{main_thm}
Let $a_1,\ldots,a_d \in \F \setminus \{0\}$ for $d \ge 4$ be such that $a_1+\ldots +a_d = 0$. There exists $ c_q < q$ such that  for $n$ sufficiently large,  any set $A \subset \F^n$ with $|A| > c_q^n$ contains \textbf{all distinct}
$x_1,\ldots,x_d$  satisfying 
\begin{equation}
\label{eq1}
a_1x_1 + a_2 x_2 + \ldots + a_d x_d = 0.
\end{equation}
\end{thm}
\medskip

\begin{ack}
The authors thank the Australian Mathematical Sciences Institute (AMSI) for their support of the second author through the Vacation Research Scholarship Program.  The first author was supported by the Australian Research Council grant DP210100162.
\end{ack}

\section{Slice rank of $d$-tensors versus Tao slice rank of functions on $A^d$}

Let $A \subset \F^n$, and $F:A^d \to \F$ be a function for $d \ge 2$. Following Tao in \cite{Tao}, we define the Tao slice rank of $F$, denoted as $\text{T-srank}(F)$, as follows: 
\medskip

If there exist functions $f:A \to \F$, $g:A^{d-1} \to \F$ and $i \in \{1,\ldots,d\}$ such that $F(x_1,\ldots,x_d) = f(x_i)g(\mathbf{x}^i)$, where $\mathbf{x}^i$ is the $(d-1)$-tuple obtained from $(x_1,\ldots,x_d)$ by removing $x_i$, then we define that $\text{T-srank}(F) = 1$. We define that $\text{T-srank}(F) \le k$ for $k \ge 2$ if there exist functions $F_1,\ldots,F_k:A^d \to \F$  of Tao slice rank equal to $1$, and such that $F = \sum_{j=1}^k F_j$.
\medskip

Denote by $V$ a finite dimensional vector space over $\F$. Any multilinear function $T:V^d \to \F$ is called a $d$-tensor. The slice rank of a $d$-tensor $T$, denoted as $\text{s-rank}(T)$, is defined as follows:
\medskip

We define that $\text{s-rank(T)} = 1$ if there exist $1$-tensor $T_1:V \to \F$ and $(d-1)$-tensor $T_2:V^{d-1} \to \F$, and $i \in \{1,\ldots,d\}$ such that $T(v^1,\ldots,v^d) = T_1(v^i) T_2(\mathbf{v}^i)$, where $\mathbf{v}^i \in V^{d-1}$ obtained by removing $v^i$ from $(v^1,\ldots,v^d)$, for any $(v^1,\ldots,v^d) \in V^d$. We define that $\text{s-rank}(T) \le k$ for $k \ge 2$ if there exist $k$ $d$-tensors $T^1,\ldots,T^k:V^d \to \F$ of the slice rank equal to $1$ such that $T = \sum_{i=1}^{k} T^i$. 
\medskip

For every $F:A^d \to \F$ we correspond the $d$-tensor $T_F$ on the space $V $ of all functions from $A$ to $\F$ defined as follows:
\begin{equation}
\label{tensor}
T_F(f_1,\ldots,f_d) := \sum_{(x_1,\ldots,x_d) \in A^d} F(x_1,\ldots,x_d) f_1(x_1)\ldots f_d(x_d).
\end{equation}
The next lemma follows immediately from the definitions of the slice rank of $d$-tensors, the Tao slice rank of functions on $A^d$ and the relationship between $F$ and $T_F$.
\begin{lem}\label{ranks_lemma}
\[
\textnormal{s-rank}(T_F) \le \textnormal{T-srank}(F).
\]
\end{lem}

Assume that  $T$ is a $d$-tensor on the space $V = \F^N$.  By multilinearity of $T$ we have:
\[
T(x^1,\ldots,x^d) = \sum_{\alpha \in [N]^d} c_{\alpha} x^1_{\alpha_1}\ldots x^d_{\alpha_d}, \mbox{    for  } (x^1,\ldots,x^d) \in V^d,
\]
where any vector $v \in V$ is represented in the coordinates as $v= (v_1,\ldots,v_N)$, and $\alpha \in [N]^d = \{1,\ldots,N\}^d$ has coordinates $\alpha_1,\ldots,\alpha_d$.
We define that a set $\mathcal{I} \subset \{1,\ldots,N\}$ is an independent set for $T$ if for any $\alpha = (\alpha_1,\ldots,\alpha_d) \in \mathcal{I}^d$ such that $c_{\alpha} \neq 0$ we have that $\alpha_1 = \ldots = \alpha_d$. 

\begin{thm}[Lovett \cite{Lov}, Theorem 1.7]
\label{indep_set_thm}
There exists a constant $c = C(d,q)$ such that for any $d$-tensor $T$ we have 
\[
\textnormal{s-rank}(T) \ge c | \J |,
\]
for any independent set $\J \subset \{1,\ldots,N\}$.
\end{thm}
\medskip

\section{Proof of Theorem \ref{cap_set_thm}}

Let $\delta > 0$ and $\mathbf{a} = (a,b,c) \in (\F \setminus \{0 \})^3$ be such that $a+b+c = 0$. Let $\epsilon > 0$, and we assume that $A \subset \F^n$ is an $(\epsilon,\delta)-$cap set. I.e.,  there exists a subset $A' \subset A$ with $|A'| \ge \delta |A|$ and such that for every $x \in A'$ there are at most $|A|^{\epsilon}$ pairs $(y,z) \in A^2$ with 
\[
ax+by+cz= 0.
\]
Denote by $F(x,y,z) = \delta_{\mathbf{o}^n}(ax+by+cz)$. 
Then we have $F(x,y,z) = \sum_{\alpha \in A^3} c_{\alpha} \delta_{\alpha_1}(x) \delta_{\alpha_2} (y) \delta_{\alpha_3}(z)$, for 
\[
c_{\alpha} = \begin{cases}
1 &\mbox{if } a \alpha_1 + b \alpha_2 + c \alpha_3 = 0 \\
0& \mbox{otherwise}. 
\end{cases}
\]
Following the construction (\ref{tensor}), the corresponding $3$-tensor $T_F$ on  $V = \{f: A \to \F \}$ is equal to 
\[
T_F(f_1,f_2,f_3) = \sum_{\alpha \in A^3} c_{\alpha} f_1(\alpha_1) f_2(\alpha_2) f_3(\alpha_3).
\]
By the assumption on $A' \subset A$, we deduce that 
\begin{equation}
\label{upper_bound1}
|\{ c_{\alpha} \neq 0 \, | \, \alpha \subset (A')^3 \}| \le \delta^{-1} |A'|^{1+\eps}.
\end{equation}
 Using Caro-Wei lower bound \cite{Car},\cite{Wei} on the independence number in a $3$-uniform hypergraph, there exists $\I \subset A'$ an independent set satisfying 
 \[
 | \I | \ge C_1 \sum_{x \in A'} \frac{1}{(d_x + 1)^\frac{1}{3}},
 \]
 where  $d_x = |\{ c_{\alpha} \neq 0 \, | \, \alpha_1 = x , \alpha \subset (A')^3\}| $ and $C_1> 0$ is a constant. The right hand side is minimised whenever all the summands are equal. Also, it follows from (\ref{upper_bound1}) that $\sum_{x \in A'} d_x \le  \delta^{-1} |A'|^{1+\eps}$.
  Therefore, there exists a constant $C_2 > 0$ such that 
 \[
  | \I | \ge C_2 |A'|^{1-\epsilon/3}.
 \]
 Finally, using Theorem \ref{indep_set_thm} and Lemma \ref{ranks_lemma}, there exists a constant $C_3 = C(q,d,\delta) > 0$ such that 
 \[
 \text{s-rank}(F) \ge C_3 |A|^{1-\epsilon/3}.
 \]
 On the other hand, using the fact that $F(x,y,z)$ is a polynomial in the coordinates of $x,y$ and $z$, it was proved in \cite{EG} that $\text{s-rank}(F) < b_q^n$ for a positive constant $b_q < q$. Finally, we  choose $\epsilon > 0$ to satisfy
 \[
 b_q^{\frac{1}{1 - \epsilon/3}} < q,
 \]
 and take any $c_q$ that satisfies $ b_q^{\frac{1}{1 - \epsilon/3}} < c_q < q$. 
Then the statement of the Theorem holds true for chosen $\epsilon$ and $c_q$.

\qed

\section{Proof of Theorem \ref{main_thm}}
By rearranging, if necessary, we always can assume that for every $k \le d-2$ we have that $a_1+\ldots +a_k \neq 0$. Let us denote by $b_{k}$, $k = 2,\ldots,d-2$, the quantities 
\[
b_k = a_1+\ldots+a_k.
\]
We apply Theorem \ref{cap_set_thm} iteratively on the equations:
\begin{equation}
\label{eq_{d-2}}
b_{d-2}t_{d-2} + a_{d-1} x_{d-1} + a_d x_d = 0,
\end{equation}
and
\[
\begin{matrix}
b_{d-3} t_{d-3} + a_{d-2} x_{d-2} = b_{d-2}t_{d-2}\\
\vdots\\
b_k t_k + a_{k+1} x_{k+1} = b_{k+1}t_{k+1}\\
\vdots\\
b_{2} t_{2} + a_{3} x_{3} = b_{3}t_{3}\\
\end{matrix}
\]
and
\[
a_1 x_1 + a_2 x_2 = b_{2}t_2.
\]
Fix $\delta > 0$. By Theorem \ref{cap_set_thm} there exists $\epsilon > 0$, such that for sufficiently large  $n$  there exists a set $A_{d-2}\subset A$ with $|A_{d-2}| \ge (1-\delta) |A|$ and such that for every $t_{d-2} \in A_{d-2}$ there are at least $|A|^{\epsilon}$ pairs $(x_{d-1},x_{d}) \in A^2$ satisfying the equation (\ref{eq_{d-2}}). Applying Theorem \ref{cap_set_thm} once again, there exists $A_{d-3} \subset A_{d-2}$ with $|A_{d-3}| \ge (1-\delta) |A_{d-2}|$ such that for every $t_{d-3} \in A_{d-3}$ there are at least $|A_{d-2}|^{\epsilon}$ pairs $(x_{d-2},t_{d-2}) \in A_{d-2}^2$ satisfying the equation 
\[
b_{d-3} t_{d-3} + a_{d-2} x_{d-2} = b_{d-2} t_{d-2}.
\]
In such way we construct a chain of subsets $A_{d-2}\supset A_{d-3} \supset \ldots \supset A_{2}$ with $|A_{k-1}| \ge (1-\delta)|A_k|$, for $k=2,\ldots,d-2$. Each $A_k$ satisfies that for any  $t_k \in A_k$ there exist at least $|A_{k+1}|^{\epsilon}$ pairs $(x_{k+1},t_{k+1}) \in A_{k+1}^2$ satisfying the equation
\[
b_{k} t_k + a_{k+1} x_{k+1} = b_{k+1} t_{k+1}.
\]
Using Theorem \ref{cap_set_thm} again, there exists $A_1 \subset A_2$ with $|A_1| \ge (1-\delta) |A_2|$ and such that for every $t_2 \in A_1$ there exist at least $|A_2|^{\epsilon}$ pairs $(x_1,x_2) \in A_2^2$ satisfying 
\[
a_1 x_1 + a_2 x_2 = b_2 t_2.
\]
Finally, we construct a solution for (\ref{eq1}) consisting of distinct elements of $A$ as follows. 
Take a pair of distinct $(x_1,x_2) \in A_2^2 \subset A^2$ satisfying that $a_1 x_1 + a_2 x_2 = b_2 t_2$ for $t_2 \in A_2$. Then there exist at least $|A_3|^{\epsilon}$ pairs $(x_3,t_3) \in A_3^2 \subset A^2$ satisfying $b_2t_2+a_3x_3 = b_3 t_3$ for $t_2$ that we already chosen. Find $(x_3,t_3)$ among these solutions such that $x_3 \not \in \{x_1,x_2\}$. Assume that we already constructed distinct $\{x_1,\ldots,x_k\} \in A^k$ satisfying that $a_1x_1 + \ldots +a_k x_k = b_k t_k$, for some $t_k \in A_k$. Since there exist at least $|A_{k+1}|^{\epsilon}$ pairs $(x_{k+1},t_{k+1}) \in A_{k+1}^2$ satisfying 
\[
b_k t_k + a_{k+1} x_{k+1} = b_{k+1} t_{k+1},
\]
we can choose one of the solutions $(x_{k+1},t_{k+1}) \in A_{k+1}^2$ satisfying that $x_{k+1} \not \in \{x_1,\ldots,x_k\}$. Notice that there exists $t_{k+1} \in A_{k+1}$ such that the sequence $(x_1,\ldots,x_{k+1}) \in A^{k+1}$ satisfies
\[
a_1x_1+\ldots + a_{k+1}x_{k+1} = b_{k+1} t_{k+1}.
\]
We  continue this process till we reach distinct $\{x_1,\ldots,x_{d-2}\} \in A$ satisfying 
\[
a_1x_1 + \ldots + a_{d-2} x_{d-2} = b_{d-2} t_{d-2}
\]
for some $t_{d-2} \in A_{d-2}$. Since there are at least $|A|^{\epsilon}$ pairs $(x_{d-1},x_d) \in A^2$ satisfying 
\[
b_{d-2} t_{d-2} + a_{d-1} x_{d-1} + a_d x_d = 0,
\]
we can choose the solution $(x_{d-1},x_d) \in A^2$ such that $x_{d} \neq x_{d-1}$ and $x_{d-1}, x_d \not \in \{x_1,\ldots,x_{d-2}\}$. This finishes the proof of the Theorem.
\qed

\begin{rema}
It seems to be more natural to try to prove Theorem \ref{main_thm} using directly Tao's reformulation \cite{Tao} of Croot-Lev-Pach-Ellenberg-Gijswijt approach together with Lovett's lower bound (Theorem \ref{indep_set_thm}) on the slice rank of the corresponding tensor. Unfortunately, this does not work out, since the lower bound on the slice rank that we obtain is not sufficiently strong.
\end{rema}


\medskip


\begin{thebibliography}{99}
\bibitem{Car} Caro Y., \emph{ New Results on the Independence Number}, Technical Report, Tel-Aviv University,
1979.
\bibitem{CLP} Croot, E.,  Lev, V. F., Pach, P.,
\emph{ Progression-free sets in $\Z_4^n$ are exponentially small. }Ann. of Math. (2) 185 (2017), no. 1, 331--337.
\bibitem{EG}  Ellenberg, J.,  Gijswijt, D.,
\emph{On large subsets of $\F^n$ with no three-term arithmetic progression.  }Ann. of Math. (2) 185 (2017), no. 1, 339--343.
\bibitem{Lov} Lovett, S.,
\emph{The analytic rank of tensors and its applications. } Discrete Anal. 2019, Paper No. 7, 10 pp. 
\bibitem{Tao}
Tao, T. ,
\emph{A symmetric formulation of the Croot-Lev-Pach-Ellenberg-Gijswijt capset bound}, 2016, \url{http://terrytao.wordpre.com/2016/05/18/a}.
\bibitem{Wei} Wei V. K.,  \emph{A Lower Bound on the Stability Number of a Simple Graph}, Technical memorandum,
TM 81--11217--9, Bell laboratories, 1981.
\end{thebibliography}
\end{document}